\documentclass[11pt]{article}
\usepackage{mathrsfs}
\usepackage{amsmath}
\usepackage{amssymb}
\usepackage{graphicx}
\usepackage{epic}
\usepackage{amsthm,latexsym}
\renewcommand{\paragraph}{\roman{paragraph}}
\def \G{\Gamma}
\newtheorem{theorem}{\scshape \mdseries  Theorem}[section]
\newtheorem{lemma}[theorem]{\scshape \mdseries  Lemma}

\newtheorem{defi}[theorem]{\scshape \mdseries Definition}
\DeclareMathOperator{\sgn}{sgn}
\begin{document}

\title{\vspace{10mm} \sf  The nullity of unicyclic signed graphs\thanks{Supported by 
 National Natural Science Foundation of China (11071002),
Program for New Century Excellent Talents in University, 
Key Project of Chinese Ministry of Education (210091), 
Specialized Research Fund for the Doctoral Program of Higher Education (20103401110002), 
Science and Technological Fund of Anhui Province for Outstanding Youth  (10040606Y33), 
NSF of Department of Education of Anhui Province (KJ2010B136), 
Project of Anhui Province for Excellent Young Talents in Universities (2009SQRZ017ZD),
Scientific Research Fund for Fostering Distinguished Young Scholars of Anhui University(KJJQ1001), 
Academic Innovation Team of Anhui University Project (KJTD001B), 
Fund for Youth Scientific Research of Anhui University(KJQN1003),
and Innovation Fund for Graduates of Anhui University.}}
\author{Yi-Zheng Fan\thanks{Corresponding author. 
Email addresses:   fanyz@ahu.edu.cn (Y.-Z. Fan),  wangyue0903@126.com (Y. Wang),  wangy@ahu.edu.cn (Y. Wang).}, \ Yue Wang, 
\ Yi Wang\\
 {\small  \sl School of Mathematical Sciences, Anhui University, Hefei 230039, P. R. China} 
}
\date{}
\maketitle

\noindent {\bf Abstract}
In this paper we introduce the nullity of signed graphs, and give some results on the nullity of signed graphs with pendant trees.
We characterize the unicyclic signed graphs of order $n$ with nullity $n-2,n-3,n-4,n-5$ respectively.



\section{Introduction}
Let $G=(V,E)$ be a  simple graph  with vertex set $V=V(G)=\{v_1,v_2,\ldots,v_n\}$ and edge set $E=E(G)$. 
The {\it adjacency matrix} $A=A(G)=(a_{ij})_{n\times n}$ of $G$ is defined as follows:
   $a_{ij}=1$ if there exists an edge  joining $v_i$ and $v_j$, and $a_{ij}=0$ otherwise. 
The {\it nullity of a simple graph} $G$, denoted by $\eta(G)$, is the multiplicity of the eigenvalue zero in the spectrum of $A(G)$.
The graph $G$ is called {\it singular} (or {\it nonsingular}) if $A(G)$ is singular  or $\eta(G)>0$ (or $A(G)$ is nonsingular or $\eta(G)=0$). 

Recently the nullity of simple graphs has been received a lot of attention.
Collatz and Sinogowitz \cite{collatz}  posed the problem of characterizing nonsingular or singular graphs.
If $G$  is a nonsingular bipartite graph, then, as shown in \cite{long}, the alternant hydrocarbon corresponding to $G$ is unstable. 
The problem is also of interest in mathematics itself, as it is closely related to the minimum rank problem of symmetric
matrices whose patterns are described by graphs \cite{fallat}.

It is known that $0\le \eta(G)\le n-2$ if $G$ is a simple graph of order $n$  containing at least one edge.
Cheng and Liu \cite{chengl} characterize the graphs of order $n$ with nullity $n-2$ or $n-3$.
Cheng, Huang and Yeh \cite{chenghy} characterize the graphs of order $n$ with nullity $n-4$.
The characterization of graphs of order $n$ with nullity $n-5$ or more is still open.
Much work is devoted to the nullity of special classes of graphs; see \cite{fan, fiorini, gong, guo, gut, hu, lics, li, lic, nath, tan, zhu}.

In this paper we discuss the nullity of the signed graphs.
A {\it signed graph} is a graph with a sign attached to each of its edges.
Formally, a signed graph $\Gamma=(G,\sigma)$ consists of a simple graph $G=(V,E)$,
  referred to as its underlying graph, and a mapping $\sigma: E \rightarrow \{+,-\}$, the edge labelling.
To avoid confusion, we also write $V(\G)$ or $V(G)$ instead of $V$, $E(G)$ instead of $E$, and $E(\G)=E^\sigma$.
The {\it adjacency matrix} of the signed graph $\G$ is  $A(\Gamma)=(a_{ij}^\sigma)$ with
$a_{ij}^\sigma=\sigma(v_iv_j)a_{ij}$, where $(a_{ij})$ is the adjacency matrix of the underlying graph $G$.
In the case of $\sigma=+$ being an all-positive edge labelling, 
   then the adjacency matrix $A(G,+)$ is exactly the classical adjacency matrix of $G$.
The {\it nullity of a signed graph} $\G$ is defined as the multiplicity of the eigenvalue zero in the spectrum of $A(\G)$, 
  and is still denoted  by $\eta(\G)$.

Signed graphs were introduced by Harary \cite{har} 
   in connection with the study of the theory of social balance in social psychology (see \cite{der}).
The matroids of graphs were extended to those of signed graphs by Zaslavsky \cite{zas}, 
and the Matrix-Tree Theorem for signed graphs was obtained by Zaslavsky \cite{zas} and by Chaiken \cite{cha}. 
More recent results on signed graphs can be found in \cite{cam,hou,hou2}.

In chemical signed graph theory, the edge  signed graph (exactly the signed graph defined here) was introduced according to
the eigenvectors (or molecular orbitals) of the adjacency matrix of the underlying graph. 
The {\it net sign} is defined to the sum of all signs of the edges of the signed graph, 
which is used to reflect the bonding capacity and rationalize the scheme of chemical bonding; see \cite{gutlee,leeli, leelu,salee} for details.

Let $\Gamma=(G,\sigma)$ be a signed graph.
The {\it sign} of a cycle $C$ of $\Gamma$ is denoted and defined by $\sgn(C)=\Pi_{e\in C} \sigma(e)$. 
The cycle $C$ is said {\it positive} or {\it negative} if $\sgn(C)=+$ or $\sgn(C)=-$.
A signed graph is said to be {\it balanced} if all its cycles are positive, 
  or equivalently, all cycles have even number of negative edges; otherwise it is called {\it unbalanced}.
There have been a variety of applications of balance; see \cite{rob}.

Suppose $\theta: V(G) \to \{+,-\}$ is any sign function.
Switching $\G$ by $\theta$ means forming a new signed graph $\G^\theta=(G, \sigma^\theta)$ whose underlying
graph is the same as $G$, but whose sign function is defined on an edge $uv$ by
$\sigma^\theta(uv)=\theta(u)\sigma(uv)\theta(v)$.
Note that switching does not change the signs or balance of the cycles of $\G$.  
If we define a (diagonal) signature matrix $D^\theta$ with $d_v=\theta(v)$ for each $v \in V(G)$, then
$A(\G^\theta)=D^\theta A(\G)D^\theta$.
Two graphs $\G_1,\G_2$ are called {\it switching equivalent}, denoted by $\G_1 \sim \G_2$, if there exists a switching function $\theta$ such that
  $\G_2=\G_1^\theta$, or equivalently $A(\G_2)=D^\theta A(\G_1)D^\theta$.

\begin{theorem} {\em \cite{hou}} \label{balance}
Let $\G$ be a signed graph. 
Then $\Gamma$ is balanced if and only if $\G=(G, \sigma) \sim (G,+)$.
\end{theorem}

Note that switching equivalence is a relation of equivalence, and two switching equivalent graphs have same nullities.
So, when we discuss the nullity of signed graphs, we can choose an arbitrary representative of each switching equivalent class.
For the unicyclic graphs, there are exactly two  switching equivalent classes.
If a unicyclic signed graph is balanced, by Theorem \ref{balance}, it is switching equivalent to one with all edges positive.
Otherwise, it is switching equivalent to one with  exactly one (arbitrary) negative edge on the cycle, by the following lemma.

 \begin{lemma} \label{unbcyc}
Let $\G$ be  an unbalanced signed unicyclic graph of order $n$.
Then $\G$ is switch equivalent to a signed unicyclic graph of order $n$ with exactly one (arbitrary) negative edge on the cycle.
\end{lemma}

{\bf Proof:}
Let $e$ be an arbitrary edge on the cycle of $\G$.
Observe that $\G-e$ is balance. 
So, by Theorem \ref{balance}, there exists a sign function $\theta$ such that $(\G-e)^\theta$ consisting of positive edges.
Returning to the graph $\G^\theta$, the edge $e$ must have negative sign as switching does not change the sign of a cycle.
The result follows.
\hfill $\blacksquare$

In this paper we concern the nullity of unicyclic signed graphs of order $n$, and characterize the unicyclic signed graphs of order $n$
with nullity $n-2,n-3,n-4,n-5$ respectively.

\section{Preliminaries}
We first introduce some concepts and notations of signed graphs.
However these definitions are based only on the underlying graph of the signed graph.
Let $\G = (G, \sigma)$ be a signed graph.
The graph $\G$ is said {\it acyclic} (respectively, {\it unicyclic}) if it contains no cycles (respectively, contains exactly one cycle). 
Particularly, the unicyclic graphs consider here are all connected.

A vertex of $\G$ is called {\it pendant} if it has degree one, and is called {\it quasi-pendant} if it is adjacent to a pendant vertex.
An edge subset $M \subseteq E(\G)$ is called a {\it matching} of $\G$ if no two edges of $M$ share a common vertex. 
A matching $M$ is called {\it maximum} in $\G$ if it has maximum cardinality among all matchings of $\G$, 
   and is called {\it perfect} if every vertex of $\G$ is incident with (exactly) one edge in $M$. 
Obviously, a perfect matching is a maximum matching. 
The cardinality of a maximum matching is called the {\it matching number} of $\G$, denoted by $\mu(\G)$.
Denote by $\G-U$, where $U\subseteq V(\G)$, the graph obtained from $\G$ 
   by removing the vertices of $U$ together with all signed edges incident to them. 
Sometimes we use the notation $G-G_1$ instead of $G-V(G_1)$, where $G_1$ is a subgraph of $G$. 

The  union of two disjoint graphs $G_1$ and $G_2$ is denoted by $G_{1} \oplus G_{2}$. 
The union of $k$ disjoint copies of $G$ is often written as $kG$. 
Denote by $K_n,K_{1,n-1},C_n$ a complete graph,  a star and a cycle all of order $n$, respectively.

Note that for a balanced signed graph $\G=(G,\sigma)$, the matrix $A(\G)$ is similar to $A(G)$ via a signature matrix 
   by Theorem \ref{balance}.
So the nullity results for simple graphs still hold for balanced signed graphs.

\begin{lemma} {\em \cite{cve2}} \label{trenul}
If $T$ is a acyclic signed graph or a signed tree of order $n$. Then
$\eta(T)=n-2\mu(T)$.
\end{lemma}

\begin{lemma} {\em \cite{cve2}} \label{bcynul}
Let $C_n$ be a balanced cycle. 
Then $\eta(C_n)=2$ if $n \equiv 0(\!\!\!\mod 4)$ and $\eta(C_n) = 0$ otherwise.
\end{lemma}

The following result can be obtained from Proposition 2.2 of \cite{fanmix}.

\begin{lemma} \label{ubcynul}
Let $C_n$ be a unbalanced signed cycle. 
Then it has eigenvalues $ 2\cos \frac{(2k-1)\pi}{n}$, $i=1,2,\ldots,n$.
Hence, $\eta(C_n)=2$ if $n \equiv 2 (\!\!\!\mod 4)$, and $\eta(C_n)=0$ otherwise.
\end{lemma}

\begin{lemma} \label{reduce}
Let $\G$ be a signed graph containing a pendant vertex, and let $H$ be the induced subgraph
of $G$ obtained by deleting this pendant vertex together with the vertex adjacent to it. Then
$$\eta(\G)=\eta(H).$$
\end{lemma}

{\bf Proof:}
Let $u_1$ be a pendant vertex of $G$, and let $u_2$ be its neighbor.
Denote $r(A)$ the rank of the matrix $A$.
Then 
$$r(A(G))= r \left(\left[\begin{array}{ccc}
			0 & \sgn(u_1u_2) & \mathbf{0}\\
			\sgn(u_1u_2) & 0 & \alpha\\
			\mathbf{0}^{T} & \alpha^{T} & A(H)
			\end{array}\right]\right)
=r \left(\left[\begin{array}{ccc}
0 & \sgn(u_1u_2) & \mathbf{0}\\
\sgn(u_1u_2) & 0 & \mathbf{0}\\
\mathbf{0}^{T} & \mathbf{0}^{T} & A(H)
\end{array}\right]\right).$$
So $r(A(G))=2+r(A(H))$, and hence $\eta(G)=\eta(H)$.  \hfill$\blacksquare$

\vspace{3mm}
The result of Lemma \ref{reduce} for simple graphs (\cite{cve2}) has been widely used for discussion of nullity.
In \cite{gong}, the authors extend Lemma \ref{reduce} to graphs with pendant trees.
We now adopt the terminologies of $k$-joining graph and pendant tree used in \cite{gong}.

\begin{defi} Let $\G_1$ be a signed graph containing a vertex $u$, 
   and let $\G_2$ be a signed graph of order $n$ that is disjoint from $\G_1$.
For $1 \le k \le n$, the $k$-joining graph of $\G_1$ and $\G_2$ with respect to $u$, 
   denoted by $\G_1(u)\odot^k \G_2$, is obtained from $\G_1 \cup \G_2$ 
   by joining $u$ and arbitrary $k$ vertices of $G_2$ with signed edges.
\end{defi}

In above definition, if $\G_1$ is a tree,  then $G_1$ is called a {\it pendant tree} of $\G_1(u)\odot^k \G_2$, 
  and $\G_1(u)\odot^k \G_2$ is said a signed graph with pendant tree.
Before we discuss the nullity of signed graphs with pendant trees, we need some notions and lemmas used in \cite{gong}.
For a signed tree $T$ on at least two vertices, a vertex $v \in T$ is called {\it mismatched} in $T$ 
   if there exists a maximum matching $M$ of $T$ that does not cover $v$; 
   otherwise, $v$ is called {\it matched} in $T$. 
If a tree consists only one vertex, then this vertex is considered mismatched.

\begin{lemma}{\em \cite{gong}} \label{mismat}
Let $T$ be a tree containing a vertex $v$. The following are equivalent:

$(1)$  $v$ is mismatched in $T$;

$(2)$ $\mu(T-v)=\mu(T)$;

$(3)$ $\eta(T-v)=\eta(T)-1.$

\end{lemma}

\begin{lemma} {\em \cite{gong}} \label{quamat}
If $v$ is a quasi-pendant vertex of a tree $T$, then $v$ is matched in $T$.
\end{lemma}

\begin{lemma} {\em \cite{gong}} \label{misnei}
If $v$ is a mismatched vertex of a tree $T$, then for any neighbor $u$ of $v$, 
  $u$ is matched in $T$, and is also matched in the component of $T-v$ that contains $u$.
\end{lemma}

The following two theorems for simple graphs were given in \cite{gong}. Here we extend them to signed graphs with a very similar proof.
\begin{theorem} \label{reduceM}
Let $T$ be a signed tree with a matched vertex $u$ and let $\G$ be a signed graph of order $n$. Then
for each integer $k \;(1\leq k\leq n)$,
$$\eta(T(u)\odot^{k}\G)=\eta(T)+\eta(\G).$$
\end{theorem}

{\bf Proof:}  
We prove the result by induction on the matching number $\mu(T)$. 
Note that $\mu(T)\geq1$ as $u$ is matched in $T$. 
If $\mu(T) = 1$, then, by Lemma \ref{trenul}, $T$ contains $p+2$ vertices and $T$ is the star $K_{1,p+1}$, where $p =\eta(T)$.
Therefore, $u$ is the unique quasi-pendant vertex of $T$. 
Suppose $v$ is a pendant vertex of $T$ that is adjacent to $u$. 
Then by Lemma \ref{reduce}, 
$$\eta(T(u)\odot^{k}\G)=\eta((T(u)\odot^{k}\G)-v-u)=\eta(pK_{1}+\G)=p+\eta(\G)=\eta(T)+\eta(\G).$$

Suppose the assertion holds for any tree $T$ with $\mu(T) \leq t \;(t \geq 1)$. 
Now we consider a tree $T$ with $\mu(T)=t+1\geq2$. 
As $\mu(T)\geq2$, we may assume that $T$ contains a pendant vertex $v$ and a quasi-pendant vertex $w$ adjacent to $v$,
  where $v,w$ are both different to $u$. 
Let $T_1 = T-v-w$. 
Then $\mu(T_{1})=t$ and $\eta(T_1)=\eta(T)$ by Lemma \ref{reduce}. 
In addition, $u$ is also matched in $T_1$. 
By Lemma \ref{reduce} and by induction, we
have
$$\eta(T(u)\odot^{k}\G)=\eta(T(u)\odot^{k}\G-v-w)=\eta(T_1(u)\odot^{k}\G)=\eta(T_1)+\eta(G)=\eta(T)+\eta(G).$$
The result follows. \hfill$\blacksquare$

\begin{theorem} \label{reduceUM}
Let $T$ be a signed tree with a mismatched vertex $u$ and let $\G$ be a signed graph of order $n$.
Then for each integer $k \;(1\leq k\leq n)$,
$$\eta(T(u)\odot^{k}\G)=\eta(T-u)+\eta(\G+u)=\eta(T)+\eta(\G+u)-1,$$
where $\G+u$ is the subgraph of $T(u)\odot^{k}\G$ induced by the vertices of $\G$ and $u$.
\end{theorem}

{\bf Proof:} 
In the tree $T$, assume that $u_1,u_2,\ldots,u_m \;(m \geq 1)$ are all neighbors of $u$, 
  and $T_1,T_2,\ldots,T_m$ are the components of $T-u$ that contain the vertices $u_1,u_2,\cdots,u_m$ respectively. 
  By Lemma \ref{misnei}, every vertex $u_{i}$ is matched in $T_{i}$ for $i=1,2,\ldots,m$.
Then
\begin{eqnarray*}
T(u)\odot^k \G & = & T_1(u_1)\odot^1 (T(u)\odot^k \G-T_1)\\ &
= & T_1(u_1)\odot^1 \left[T_2(u_2)\odot^1 (T(u)\odot^k \G-\oplus_{i=1}^2 T_i)\right]\\
                  &=& \cdots \\
                    & = & T_1(u_1)\odot^1\left[T_2(u_2)
                    \odot^1 \cdots \odot^1 \left[T_m(u_m)\odot^1(T(u)\odot^k \G-\oplus_{i=1}^mT_i)\right]\right]\\
                  & = & T_1(u_1)\odot^1\left[T_2(u_2)
                    \odot^1 \cdots \odot^1 \left[T_m(u_m)\odot^1(\G+u)\right]\right].
                    \end{eqnarray*}
Applying Theorem \ref{reduceM} repeatedly, we have
\begin{eqnarray*}\eta(T(u)\odot^k G) & = & \eta \left(T_1(u_1)\odot^1\left[T_2(u_2)
                    \odot^1 \cdots \odot^1 \left[T_m(u_m)\odot^1(\G+u)\right]\right]\right)\\
                  &=& \eta(T_1)+\eta \left(T_2(u_2)
                    \odot^1 \cdots \odot^1 \left[T_m(u_m)\odot^1(\G+u)\right]\right) \\
		    &=& \cdots \\
                  &=& \sum_{i=1}^{m-1}\eta(T_i)+\eta \left(T_m(u_m)\odot^1(\G+u)\right)\\
                  &=& \sum_{i=1}^{m}\eta(T_i)+\eta(\G+u).
                  \end{eqnarray*}
As $u$ is mismatched in $T$, by Lemma \ref{mismat},
$\sum_{i=1}^{m}\eta(T_i)=\eta(T-u)=\eta(T)-1$. \hfill $\blacksquare$

\section{Nullity of unicyclic signed graphs}

Let $\G$ be a unicyclic signed graph and let $C$ be the unique cycle of $\G$. 
 For each vertex $v\in C$, denote by $\G\{v\}$ an induced connected subgraph of $\G$ with maximum possible of vertices,
    which contains the vertex $v$ and contains no other vertices of $C$.
One can find that $\G\{v\}$ is a tree and $\G$ is obtained by identifying the vertex $v$ of $\G\{v\}$ with the vertex $v$ on $C$
  for all vertices $v\in C$. 
The unicyclic signed graph $\G$ is said of {\it Type I} if there exists a vertex $v$ on the cycle such that $v$ is
    matched in $\G\{v\}$; otherwise, $\G$ is said of {\it Type II}.

If $\G$ is of Type I, then $\G=\G\{v\}\odot^{2}(\G-\G\{v\})$ for some matched vertex $v$ of $\G\{v\}$, 
   where $\G\{v\}(v)$ and $\G-\G\{v\}$ are both nontrivial graphs. 
Thus, by Theorem \ref{reduceM}, $$\eta(\G)=\eta(\G\{v\})+\eta(\G-\G\{v\}).$$

If $\G$ is of Type II and $\G$ is not a cycle, then by Lemma \ref{misnei}, 
 for each vertex $v$ on the cycle such that $\G\{v\}$ is nontrivial, every neighbor of $v$ in $\G\{v\}(v)$ is matched in
 the component of $\G\{v\}(v)-v$ that contains the neighbor. 
Note that $G\{v\}(v)-v$ may be a forest but each component of the forest contains at least two vertices by Lemma \ref{quamat}. 
By Theorem \ref{reduceUM},
$$\eta(\G)=\eta(\G\{v\}-v)+\eta((\G-\G\{v\})+v).$$
Applying Theorem \ref{reduceUM} repeatedly, we have
$$\eta(\G)=\sum_{v \in V(C)}\eta(\G\{v\}-v)+\eta(C)=\eta(\G-C)+\eta(C).$$

By the above discussion, we get the following result immediately.

\begin{theorem} \label{reduceUNI}
Let $\G$ be a unicyclic signed graph and let $C$ be the cycle of $\G$. 
If $\G$ is of Type I and let $v \in C$ be matched in $\G\{v\}$, then
$$\eta(\G)=\eta(\G\{v\})+\eta(\G-\G\{v\}).$$
If $\G$ is Type II, then
$$\eta(\G)=\eta(\G-C)+\eta(C).$$
\end{theorem}

\begin{theorem}
Let $\G$ be a  unicyclic signed graph of order $n$.Then

{\em (1)} $\eta(\G)=n-2$ if and only if $\G$ is the balanced cycle $C_4$.

{\em (2)} $\eta(\G)=n-3$ if and only if $\G$ is the cycle $C_{3}$.

\end{theorem}

{\bf Proof:} 
The sufficiency for (1) or (2) can be verified by Lemmas \ref{bcynul} and \ref{ubcynul}.
Suppose $\eta(\G)=n-2$. 
If $\G$ is exactly a cycle $C_n$, by Lemmas \ref{bcynul} and \ref{ubcynul}, $\G$ is the balanced cycle $C_4$.
Now assume $\G$ contains pendant edges and let $C_l$ be a cycle of $\G$.
If $\G$ is of type I, then  for some some vertex $v$ of the cycle matched in $\G\{v\}$,
$\G=\G\{v\}\odot^{2}(\G-\G\{v\})$, 
   where $\G\{v\}$ and $\G-\G\{v\}$ are both nontrivial trees of order $n_1$ and $n-n_1$ respectively.
By Theorem \ref{reduceUNI} and Lemma \ref{trenul},
$$\eta(\G)=\eta(\G\{v\})+\eta(\G-\G\{v\})=n_1-2 \mu(\G\{v\})+n-n_1-2(\G-\G\{v\}),\eqno(3.1)$$
which implies $\eta(\G) \le n-4$, a contradiction.
If $\G$ of type II, by the discussion prior to Theorem \ref{reduceUNI}, each component of $G-C_l$ is nontrivial.
By Theorem \ref{reduceUNI} and Lemma \ref{trenul},
$$ \eta(\G)=\eta(\G-C_l)+\eta(C_l)=n-l-2\mu(\G-C_l)+\eta(C_l), \eqno(3.2)$$
which implies $\eta(C_l)=l+2[\mu(\G-C_l)-1] \ge l \ge 3$, a contradiction.
The first claim follows.

For the second claim, if $\eta(\G)=n-3$ and $\G$ is a cycle, by Lemmas \ref{bcynul} and \ref{ubcynul}, 
 $\G$ is the cycle $C_3$.
Assume that $\eta(\G)=n-3$ and $\G$ is not a cycle.
By (3.1), $\G$ cannot be of type I .
If $\G$ is of type II, by (3.2), $\eta(C_l)=l+2\mu(\G-C_l)-3 \ge l-1 \ge 2$, with equality only if $l=3$. 
But $\eta(C_3)=0$, so this case cannot occur.
\hfill$\blacksquare$

Before we characterize the unicyclic signed graphs of order $n$ with nullity $n-4$ or $n-5$, we need to introduce 
four graphs  in Fig. 3.1, where $U_1(r,s)$ (respectively, $U_2(r,s)$), $r \ge s \ge 0, r+s \ge 1$,
  is obtained from a triangle (respectively,  a square) by attaching $r,s$ pendant edges at two vertices 
  (respectively,  two non-adjacent vertices), and $U_3(r)$ (respectively, $U_4(r)$),  $r \ge 1$,
  is obtained from  a square (respectively, a triangle) by identifying one vertex with a pendant
vertex of a star $K_{1,r+1}$.

\vspace{3mm}
\begin{center}
\includegraphics[scale=.7]{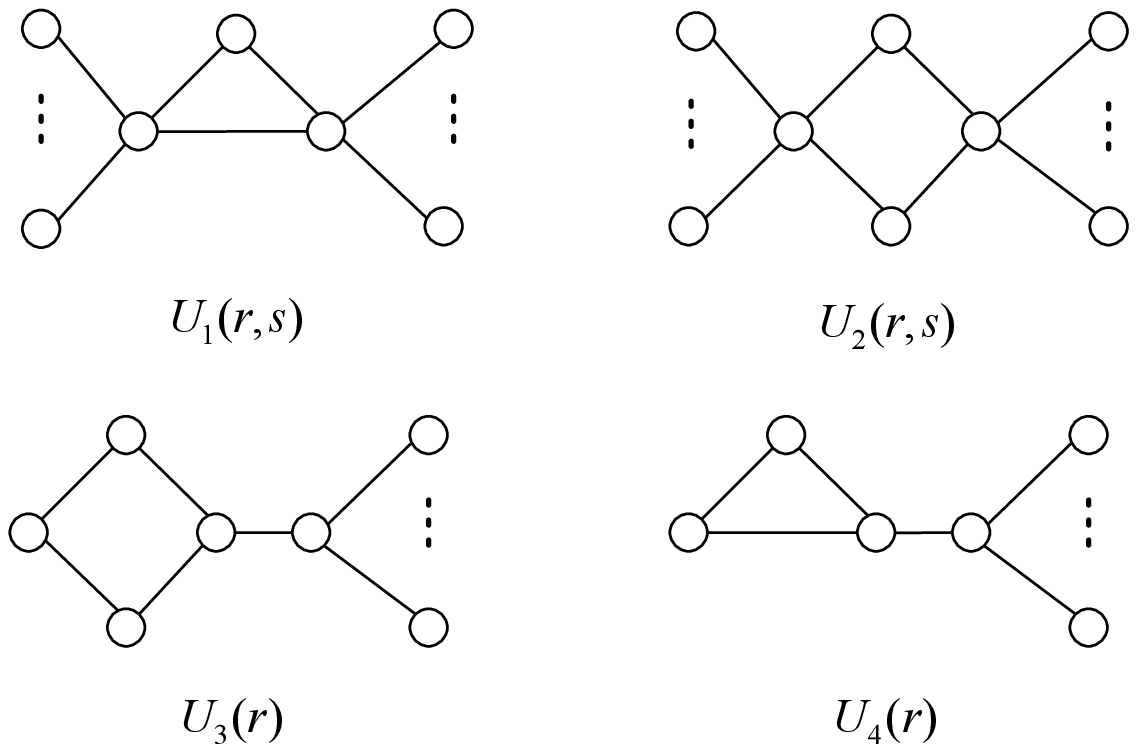}

\vspace{3mm}

\small Fig. 3.1. Four unicyclic graphs $U_1(r,s),U_2(r,s),U_3(r),U_4(s)$
\end{center}

\begin{theorem}
Let $\G$ be a  unicyclic signed graph of order $n \ge 4$. 
Then $\eta(G) = n-4$ if and only if $G$ is  one of
the following unicyclic signed graphs of order $n$: 
unbalanced $C_4$, unbalanced $C_6$,  and the signed graphs with $U_1(r,s)$ or $U_2(r,s)$ in Fig. 3.1 as underlying graph, 
the balanced signed graph with $U_3(r)$ in Fig. 3.1 as underlying graph.
\end{theorem}

{\bf Proof:}
The sufficiency is easily verified by Lemma \ref{ubcynul} and Theorem \ref{reduceUNI}.
Now suppose $\eta(\G) = n-4$.
If $\G$ is a cycle, by Lemmas \ref{bcynul} and \ref{ubcynul}, $G$ is the unbalanced cycle $C_4$ or $C_6$.
Next assume $\G$ contains pendant edges, and let $C_l$ be a cycle of $\G$.

If $\G$ is of type I, for some some vertex $v$ of the cycle matched in $\G\{v\}$,
by (3.1), $\mu(\G\{v\})=\mu(\G-\G\{v\})=1$.
So both $\G\{v\}$ and $\G-\G\{v\}$ are stars, 
   and $\G$ is obtained by identifying the center of $\G\{v\}$ with two vertices of $\G-\G\{v\}$. 
Thus $\G$ is the signed graph of order $n$ with $U_1(r,s)$ or $U_2(r,s)$ as underlying graph; see Fig. 3.1.

If $\G$ is of type II, by (3.2), $\eta(C_l)=l+2\mu(\G-C_l)-4 \ge l-2 \ge 1$.
So $\eta(C_l)=2$ by Lemmas \ref{bcynul} and \ref{ubcynul}.
So $l=4$ and $C_4$ is balance.
We also find $\mu(\G-C_4)=1$ and hence $\G-C_4$ is a star.
So $G$ is the signed graph of order $n$ with $U_3(r)$ in Fig. 3.1 as underlying graph, which is 
  obtained by identifying a a vertex of $C_4$ with one pendant vertex of a star of order $n-3$ . 
\hfill$\blacksquare$

\begin{theorem}
Let $\G$ be a  unicyclic signed graph of order $n \ge 5$. 
Then $\eta(U) = n-5$ if and only if $G$ is  one of the following graphs of order $n$: the cycle $C_5$, and signed graph with
$U_4(r)$ in Fig. 3.1 as underlying graph.
\end{theorem}

{\bf Proof:} 
The sufficiency is easily verified by Lemma \ref{ubcynul} and Theorem \ref{reduceUNI}.

Now suppose $\eta(\G) = n-5$.
If $\G$ is a cycle, by Lemmas \ref{bcynul} and \ref{ubcynul}, $G$ is the cycle $C_5$.
Next assume $\G$ contains pendant edges, and let $C_l$ be a cycle of $\G$.
If $\G$ is of type I, for some some vertex $v$ of the cycle matched in $\G\{v\}$,
by (3.1), $2[\mu(\G\{v\})+\mu(\G-\G\{v\})]=5$, impossible.

If $\G$ is of type II, by (3.2), $\eta(C_l)=l+2\mu(\G-C_l)-5 \ge l-3 \ge 0$.
If $\eta(C_l)=0$,  then $l=3$, and $\mu(\G-C_l)=1$ which implies $\G-C_l$ is a star.
So $\G$ is obtained by identifying one vertex of $C_3$ with an pendant vertex of a star of order $n-2$; see the graph $U_4(r)$ in Fig. 3.1.
If $\eta(C_l)=2$, then $l \le 5$ from the above inequalities.
By Lemmas \ref{bcynul} and \ref{ubcynul}, this case cannot occur. 
\hfill$\blacksquare$

\end{document}